\documentclass[12pt]{article}
 \usepackage{amssymb}
 \usepackage{amsmath}
 \usepackage[latin1]{inputenc}
 \usepackage[T1]{fontenc}
 \usepackage{times}
 \usepackage{hyperref} 
\date{}

 \usepackage{graphicx}
 \usepackage{epsfig}

\newtheorem{theorem}{Theorem}[section]
\newtheorem{lemma}[theorem]{Lemma}

\newtheorem{e-proposition}[theorem]{Proposition}

\newtheorem{e-definition}[theorem]{Definition\rm}
\newtheorem{remark}{\it Remark.\/}

\newtheorem{theoreme}{Th\'eor\`eme}[section]

\newtheorem{proposition}[theoreme]{Proposition}

 \newcommand{\R}{\mathbb{R}}
\newcommand{\E}{\mathbb{E}}
\begin{document}
\begin{center}
{\LARGE { Strong uniform convergence rates  of the linear wavelet estimator of a multivariate copula  density}}
 \end{center}
 \bigskip

\begin{center}
\large{Cheikh Tidiane Seck, Salha Mamane}\vspace{3mm}\\
 {\small \it D\'epartement de math\'ematiques, Universit\'e Alioune Diop, Bambey, S\'en\'egal}\\
{\small \it School of Statistics and Actuarial Science, University of the Witwatersrand, Johannesburg, South-Africa}
\end{center}
 \bigskip
\textbf{Abstract}\\
In this paper, we investigate the  almost sure convergence, in supremum norm, of the rank-based linear wavelet estimator for the multivariate copula density. Based on empirical process tools, we prove a uniform limit law for the deviation,  from its expectation, of an oracle estimator (obtained for known margins), from which we derive  the exact  convergence rate for the rank-based linear estimator. This rate reveals to be optimal in a minimax sense over Besov balls for the supremum norm loss, whenever the resolution level is suitably chosen.\\


\\

\section{Introduction}
\label{intro}
A copula  is a multivariate distribution function $C$ defined on  $[0,1]^d, d\geq 2$, with uniform margins.  Unlike the linear correlation coefficient, it gives  a full characterization of the dependence between random variables, be it linear or nonlinear.  Given a vector $(X_1,\ldots,X_d)$ 
of continuous random variables with  marginal distribution functions $F_1, \cdots, F_d$,  the copula $C$ may be defined as  the 
joint cumulative distribution function of the random vector $(F_1(X_1),\ldots,F_d(X_d))$ . 
If it exists,  the copula  density is defined as the derivative, $c$, of the copula distribution function $C$ with respect to the Lebesgue measure :
$$ c(u_1\ldots,u_d)=\frac{\partial^d}{\partial u_1,\ldots \partial u_d }C(u_1,\ldots,u_d), \enskip \forall \;(u_1,\ldots,u_d)\in (0,1)^d.$$

 Nonparametric estimation of copula density  is an active reseach domain that has been  investigated by many authors. For instance, \cite{gm} and \cite{fs}  used  convolution kernel methods to construct consistent estimators for the  copula density, while \cite{ss} employed techniques based on Bernstein polynomials. A drawback of kernel methods is the existence of boundary effects due to the compact support of the copula function. 
To overcome this difficulty, some approaches have been proposed. For example  \cite{gm}   used a mirror-reflexion technique, while \cite{ch} employed a local linear kernel procedure. In the same vein \cite{ogv} proposed some improved copula kernel estimators in order to face the boundary bias. 
Recently, \cite{gcp} introduced kernel-type estimators for the copula density, based on a probit transformation method that can take care of the boundary effects.

In this paper, 
 we deal more neatly with the boundary bias problem by using  wavelet methods, which are very convenient to describe  features of functions at the edges and corners of the unit cube, because of their good localization properties. Indeed,  wavelet bases automatically handle the boundary effects by locally adapting to the properties of the curve being estimated.  The use of wavelet methods in density and regression estimation problems is surveyed in \cite{hkt}, where approximation properties of wavelets are discussed at length. For more details on  wavelet theory we refer to \cite{mey}, \cite{dau}, \cite{mal} and \cite{vid}  and references therein.

Wavelet methods have been already used in nonparametric copula density estimation. For instance, \cite{gmt} dealt with a rank-based linear wavelet estimator of the bivariate copula density and established, under certain conditions, its optimality in the minimax sense on Besov-balls for the $L_2$-norm loss, as well as on H$\ddot{o}$lder-balls for the pointwise-norm loss.  \cite{aut} extended  these results to the nonlinear thresholded estimators of multivariate copula densities. These nonlinear estimates are near  optimal (up to a logarithmic factor) for the $L_2$-norm loss, and have the advantage of being adaptive to the regularity of the copula density function. In the same spirit, \cite{gh} established  an upper bound on $L_p$-losses, $2\le p <  \infty$, for linear wavelet-based estimators of the bivariate copula density, when the latter is bounded.

Our goal in this paper, is to establish the exact almost sure convergence rate, in supremum norm loss, for the rank-based linear wavelet estimator of  the multivariate copula density.
Our methodology is largely inspired by \cite{gn}, who established  almost sure convergence rates, in supremum norm loss, for the linear wavelet estimator of a univariate density function on $\mathbb{R}$.
Here, we want to extend this result to multivariate copula densities on $(0,1)^d$. We prove that if copula density $c$ is regular enough  (i.e. $c$ belongs to a Besov space of regularity $t$,  corresponding to the H$\ddot{o}$lder space of order $t$) and the resolution level, say $j_n$, satisfies : $2^{j_n}\simeq (n/\log n)^{1/(2t+d)}$, then the rank-based linear wavelet estimator achieves the optimal minimax rate, for supremum norm loss, over Besov-balls.


The rest of the paper is organized as follows. In Section 2, we recall some facts on wavelet theory and define  the rank-based linear wavelet estimator of the multivariate copula density as in \cite{aut}. Section 3 presents the main theoretical results along with some comments. In appendix A, we recall some useful facts on empirical process theory. Appendix B contains the proof of the uniform limit law given in Proposition \ref{p1}. 


\section{Wavelet theory and Estimation procedure}
\label{sec:1}
Let $\phi$ be a father wavelet and $\psi$ its associated mother wavelet, which are both assumed compactly supported.
\cite{cdv}  proposed  orthonormal wavelet bases for $L_2([0,1])$, the space of all measurable and square integrable functions  on $[0,1]$. Precisely, for all fixed  $j_0\in\mathbb{N}$, the family  $\{\phi_{j_0,l}:l=1,\ldots,2^{j_0}\} \bigcup \{\psi_{j,l}: j\geq j_0, l=1,\ldots,2^j \} $ is an orthonormal basis for $L_2([0,1])$, where $\phi_{j,l}(u) = 2^{j/2}\phi(2^{j}u-l)$ and $\psi_{j,l}(u) = 2^{j/2}\psi(2^{j}u-l),\; \forall j,l\in\mathbb{Z}$, $u\in [0,1]$ . Using the tensorial product, one can construct 
a multivariate wavelet basis for  $L_2([0,1]^d), d\ge 2$. For  $ \textbf{k}=(k_1,\ldots,k_d)\in \mathbb{Z}^d$, define the following functions of $\textbf{u}=(u_1,\ldots,u_d)\in [0,1]^d$ :
\begin{eqnarray*}
	\phi_{j_0,\textbf{k}}(u_1\ldots,u_d) =\prod_{m=1}^{d} \phi_{j_0,k_m}(u_{m}),\\
	\psi^{\epsilon}_{j,\textbf{k}}(u_1,\ldots,u_d) = \prod_{m=1}^{d} \phi^{1-\epsilon_m}_{j,k_m}(u_{m})\psi^{\epsilon_m}_{j,k_m}(u_{m}),
\end{eqnarray*}
where $\epsilon=(\epsilon_1,\ldots,\epsilon_d)\in \mathcal{S}_d=\{0,1\}^d\setminus\{(0,\ldots,0)\}$.
Then the family 
$ \{\phi_{j_0,\textbf{k}}, \psi^{\epsilon}_{j,\textbf{h}} : j\ge j_0, \textbf{k}\in \{1,\ldots,2^{j_0}\}^d,  \textbf{h}\in \{1,\ldots,2^j \}^d , \epsilon\in\mathcal{S}_d \}$
is an orthonormal basis for $L_2([0,1]^d)$, for any fixed $j_0\in\mathbb{N}$.
Thus, assuming that the copula density $c$ belongs to $L_2([0,1]^d)$, we have the following representation :
\begin{equation}
\label{ss02}
c(\textbf{u})=\sum_{\textbf{k}\in\{1,\ldots, 2^{j_0}\}^d }\alpha_{j_0,\textbf{k}}\phi_{j_0,\textbf{k}}(\textbf{u}) + \sum_{j\geq j_0}\sum_{\textbf{k}\in\{1,\ldots, 2^{j}\}^d }\sum_{\epsilon\in\mathcal{S}_d}\beta^{\epsilon}_{j,\textbf{k}}\psi^{\epsilon}_{j,\textbf{k}}(\textbf{u}),
\end{equation}
for all $\textbf{u}\in [0,1]^d,$ where the scaling coefficients  $\alpha_{j_0,\textbf{k}}$ and wavelet coefficients $\beta^{\epsilon}_{j,\textbf{k}}$ are respectively defined as
\begin{equation*}
\alpha_{j_0,\textbf{k}}=\int_{[0,1]^d} c(\textbf{u})\phi_{j_0,\textbf{k}}(\textbf{u})d\textbf{u}\quad \text{ and } \quad\beta^{\epsilon}_{j,\textbf{k}}=\int_{[0,1]^d} c(\textbf{u})\psi^{\epsilon}_{j,\textbf{k}}(\textbf{u})d\textbf{u}.
\end{equation*}

Now, let $(\textbf{X}_{1}, \cdots, \textbf{X}_n)$ be an independent and identically distributed (i.i.d) sample of the random vector $\textbf{X}=(X_1,\ldots,X_d)$, with continuous marginal distribution functions $F_1,\ldots,F_d$, and where $\textbf{X}_i=(X_{i1},\ldots,X_{id}), i=1,\ldots,n$. 
The distribution function of the random vector $\textbf{U}_i=(F_1(X_{i1}),\ldots, F_d(X_{id}))$ is the copula $C$ and its density, if it exists, is $c$. Denoting the expectation operator by $\mathbb{E}$, the coefficients $\alpha_{j_0,\textbf{k}}$ and $\beta^{\epsilon}_{j,\textbf{k}}$ can be rewritten as follows :
$$\alpha_{j_0,\textbf{k}}=\mathbb{E}[\phi_{j_0,\textbf{k}}(\textbf{U}_i)],\qquad \enskip \beta^{\epsilon}_{j,\textbf{k}}=\mathbb{E}[\psi^{\epsilon}_{j,\textbf{k}}(\textbf{U}_i)].$$
If the margins $F_1,\ldots,F_d$ were known, natural estimators for  $\alpha_{j_0,\textbf{k}}$ and $\beta^{\epsilon}_{j,\textbf{k}}$ would be given by, 
\begin{equation}
\tilde{ \alpha}_{j_0,\textbf{k}}=\frac{1}{n}\sum_{i=1}^n \phi_{j_0,\textbf{k}}(\textbf{U}_i ) \enskip,\enskip \tilde{ \beta}^{\epsilon}_{j,\textbf{k}}=\frac{1}{n}\sum_{i=1}^n \psi^{\epsilon}_{j,\textbf{k}}(\textbf{U}_i).
\end{equation}

But, usually the marginal ditribution functions $F_1,\ldots, F_d$ are unknown ; and it is customary to replace them by their empirical counterparts $F_{1n},\ldots,F_{dn}$ (or rescaled versions thereof), with
$$  F_{jn}(x_j)=\sum_{i=1}^n\mathbb{I}(X_{ij}\le x_j),\quad j=1,\ldots, d, \quad x_j\in\mathbb{R}, $$ where $\mathbb{I}(\cdot)$ denotes the indicator function. Then, putting  $\hat{\textbf{U}}_i= (\hat{U}_{i1},\ldots,\hat{U}_{id})$, where $\hat{U}_{ij}=F_{jn}(X_{ij}), j=1,\ldots, d$ ; $ i=1,\ldots,n$, the modified empirical coefficients are
\begin{equation}
\hat{ \alpha}_{j_0,\textbf{k}}=\frac{1}{n}\sum_{i=1}^n \phi_{j_0,\textbf{k}}(\textbf{U}_i) \enskip,\enskip \hat{ \beta}^{\epsilon}_{j,\textbf{k}}=\frac{1}{n}\sum_{i=1}^n \psi^{\epsilon}_{j,\textbf{k}}(\textbf{U}_i).
\end{equation} 

Now, choosing a suitable resolution level $j_n\ge j_0$ and considering the orthogonal projection of $c$ onto the sub-space $V_{j_n}$ of the underlying multiresolution analysis on $L_2([0,1]^d)$, we obtain the rank-based linear wavelet estimator of $c$ :
\begin{equation}
\label{ss03}
\hat{c}_{j_n}(\textbf{u})=\sum_{\textbf{k}\in\{1,\ldots,2^{j_n}\}^d}\hat{\alpha}_{j_n,\textbf{k}}\phi_{j_n,\textbf{k}}(\textbf{u}),\qquad \textbf{u}\in(0,1)^d. 
\end{equation}

As remarked in \cite{gmt}, the estimator $\hat{c}_{j_n}$ is not necessarily  a density, because it can take negative values on parts of its domain and fails to integrate to one. In practice, some truncations and normalizations are necessary for its use.\\

To obtain the exact rate of convergence of the linear estimator $\hat{c}_{j_n}$, our methodology of proof follows the empirical process approach developped in \cite{em} (see also \cite{gn},\cite{gg}). In fact, we can rewrite $\hat{c}_{j_n}$ in terms of the empirical measure. For $j_n$ fixed, define the following kernel functions : 
\begin{equation}\label{ker}
\widetilde{K}(x,y)=\sum_{l=1}^{2^{j_n}}\phi(x-l)\phi(y-l), \qquad  (x,y)\in \R
\end{equation}
$$\quad\widetilde{K}_{j_n}(x,y)=2^{j_n}\widetilde{K}(2^{j_n}x, 2^{j_n}y),\qquad  (x,y)\in \R. $$
For $\textbf{x}=(x_1,\ldots,x_m)\in {\R}^d,\; \textbf{y}=(y_1,\ldots,y_m)\in {\R}^d$, set :
\begin{equation} \textbf{K}(\textbf{x},\textbf{y})=\prod_{m=1}^d\widetilde{K}(x_m,y_m),
\end{equation}

\begin{equation} \textbf{K}_{j_n}(\textbf{x},\textbf{y})=\prod_{m=1}^d\widetilde{K}_{j_n}(x_m,y_m).
\end{equation}
Then,  the linear wavelet estimator $\hat{c}_{j_n}$ can be transformed into  
\begin{equation}\label{ss04}
\hat{c}_{j_n}(\textbf{u})=\frac{1}{n}\sum_{i=1}^n \textbf{K}_{j_n}(\hat{\textbf{U}}_i,\textbf{u})=\frac{2^{dj_n}}{n}\sum_{i=1}^n \textbf{K}(2^{j_n}\hat{\textbf{U}}_i,2^{j_n}\textbf{u}).
\end{equation}

\section{Asymptotic behaviour of the estimator}
Let's introduce an auxillary estimator $\tilde{c}_{j_n}$ corresponding to the case where the marginal distribution functions $F_1,\ldots,F_d$ are known. In this situation, $(U_{i1},\ldots,U_{id})=(F_1(X_{i1}),\ldots,F_d(X_{id})), i=1,\ldots,n$ are direct observations of the copula $C$, and $\tilde{c}_{j_n}$ may be defined as
\begin{equation}
\label{ss05}
\tilde{c}_{j_n}(\textbf{u})=\sum_{\textbf{k}\in\{1,\ldots,2^{j_n}\}^d}\tilde{\alpha}_{j_n,\textbf{k}}\phi_{j_n,\textbf{k}}(\textbf{u}), 
\end{equation}
where 
\begin{equation}
\tilde{ \alpha}_{j_n,\textbf{k}}=\frac{1}{n}\sum_{i=1}^n \phi_{j_n,\textbf{k}}(F_1(X_{i1}),\ldots, F_d(X_{id}) )
\end{equation}
is an unbiased estimator of $\alpha_{j_n,\textbf{k}}$.

For all $\textbf{u}\in (0,1)^d$, we can decompose the estimation error $\hat{c}_{j_n}-c$ as
\begin{eqnarray}\label{dec}
 \hat{c}_{j_n}(\textbf{u})- c(\textbf{u})&=& [\hat{c}_{j_n}(\textbf{u})-\tilde{c}_{j_n}(\textbf{u})] + [\tilde{c}_{j_n}(\textbf{u})-\E\tilde{c}_{j_n}(\textbf{u})]+ [\E\tilde{c}_{j_n}(\textbf{u}) -c(\textbf{u})]\nonumber \\
 &=:& R_n (\textbf{u})+D_n(\textbf{u})+B_n(\textbf{u}).
\end{eqnarray} 
To obtain the almost sure convergence rate of $\hat{c}_{j_n}$ uniformly in  $\textbf{u}\in(0,1)^d$, we have to investigate the limiting behavior of each of the three above terms. We need the following hypotheses in the sequel :
\begin{itemize}
	\item[(H.1)] The father wavelet $\phi \in L^2(\mathbb{R})$ is  bounded, compactly supported and admits a bounded derivative $\phi'$.
	\item[(H.2)] There exists a bounded and compactly supported function $\Phi: \mathbb{R}\rightarrow \mathbb{R}_{+}$ such that \\$ |\widetilde{K}(x,y)|\leq \Phi(x-y)$
	and the function $\theta_{\phi}(x)=\sum_{k=1}^{2^{j_n}}|\phi(x-k)|$ is bounded.
	\item[(H.3)] The kernel $\widetilde{K}$ satisfies : for all $y\in\mathbb{R},$ 
	 $ \int_{-\infty}^{\infty}\widetilde{K}(x,y)dx=1.$ 
	\item[(H.4)] As $n\rightarrow\infty,$  the sequence $(j_n)_{n\ge 0}$ satisfies : $$  j_n\nearrow\infty,\qquad\frac{n}{j_n2^{(d+1)j_n}}\rightarrow\infty,\qquad\frac{j_n}{\log\log n}\rightarrow\infty.$$
\end{itemize}

\begin{remark} \mbox{}
	 Hypotheses (H.1), (H.2) and (H.3) are usual conditions that are satisfied by many wavelets bases, for example the Daubechies wavelets and the Haar wavelet $\phi(u)=1_{[0,1]}(u)$. The conditions in Hypothesis  (H.4) are analogous to some conditions imposed  on the bandwidth parameter in convolution-kernel estimation methods.
\end{remark}

The following proposition gives the asymptotic behavior of the second term $D_n(\textbf{u})$, corresponding to the deviation of the auxillary estimator $\tilde{c}_{j_n}$ from its expectation. In the sequel, we denote $I=(0,1)$, $\Vert c\Vert_{\infty}=\sup_{\textbf{u}\in I^d}|c(\textbf{u})|$
 and for any real function $\varphi$ defined on $\R^d, d\ge 1$, $\Vert \varphi\Vert_{\infty}=\sup_{\textbf{x}\in \R^d}|\varphi(\textbf{x})|$.
\begin{proposition}\label{p1} Suppose that  the father wavelet $\phi$ is uniformly continuous, with compact support $[0,B]$, where $B$ is a positive integer. Further, assume that the copula density $c$ is continuous and  bounded on $I^d$ and that Hypotheses  (H.1), (H.2), (H.3)  and (H.4) hold. Then, almost surely (a.s.),
	\begin{equation}\label{rate}
	\lim_{n\rightarrow\infty}r_n \sup_{\textbf{u}\in I^d}\frac{|\tilde{c}_{j_n}(\textbf{u}) - \mathbb{E}\tilde{c}_{j_n}(\textbf{u})|}{\sqrt{\int_{{\R}^d} \textbf{K}^2(\textbf{x},2^{j_n}\textbf{u})d\textbf{x}}}=\sqrt{{\Vert c\Vert_{\infty}}},
	\end{equation}
	with  \begin{equation} r_n=\sqrt{\frac{n}{(2d\log 2) j_n 2^{d j_n}}}. 
	\end{equation}
\end{proposition} 
\textbf{Proof :} It is largely inspired by \cite{gn} and is postponed to Appendix B. It will consist in establishing a lower bound  and an upper bound for the limit in \eqref{rate}, a methodology  borrowed from \cite{em} (see also \cite{gg}). \\

\begin{remark} {\rm Proposition \ref{p1} gives the exact almost sure convergence rate, in supremum norm, of the deviation $D_n(\textbf{u})$ to zero, which is of order $O(\sqrt{j_n 2^{dj_n}/n})$. In fact, by hypotheses (H.1-2-3) the quantity $\int_{{\R}^d} \textbf{K}^2(\textbf{x},2^{j_n}\textbf{u})d\textbf{x}$ can be bounded above and below :  there exist two positive constants $D_1$ and $D_2$ independent of $\textbf{u}$ and $n$ such that :
\begin{equation}\label{borne}
D_1\le \int_{{\R}^d}\textbf{K}^2(\textbf{x},2^{j_n}\textbf{u})d\textbf{x}\le D_2.
\end{equation}
This readily implies
\begin{equation}\label{dev}
\sup_{\textbf{u}\in I^d} |D_n(\textbf{u})|=O\left(\sqrt{\frac{j_n 2^{dj_n}}{n}}\right), \quad a.s.
\end{equation} }
\end{remark}

The following theorem constitutes our principal result. We need some notation before stating it. Let $N$ be a positive integer and $t=N+\alpha, 0<\alpha\le 1$. For any bounded real function $f$ defined on $I^d$ and possessing derivatives up to order $N$, set 
\begin{equation}
\Vert f\Vert_{t,\infty,\infty}= \Vert f\Vert_{\infty}+ \sum_{k=0}^{N}\sup_{u\neq v, u,v\in I^d}\frac{|f^{(k)}(u)-f^{(k)}(v)|}{|u-v|^{\alpha}}.
\end{equation} 
We say that $f$ belongs to the Besov space  of regularity $t$, $B^t_{\infty,\infty}$, if and only if $\Vert f\Vert_{t,\infty,\infty}<\infty$.\\
 \indent The following condition is also needed for the proof :\\
\textit{Condition} $1(N) :$ the father wavelet  $\phi$ is  compactly supported and admits weak derivatives up to order $N, N\in\mathbb{N}$, that are all in $\mathcal{L}^p$ for some $1\le p\le \infty$.
\begin{theorem}\label{teo1}
Suppose that the assumptions of Proposition \ref{p1} are fulfilled.  If, moreover,  $c$ belongs to $B_{\infty,\infty}^t$ and $\phi$ satisfies Condition $1(N)$, with $0<t<N+1$,
then, as $n\rightarrow\infty$, we have
\begin{equation}\label{t1}
\sup_{\textbf{u}\in I^d} |\hat{c}_{j_n}(\textbf{u})-c(\textbf{u})|=O\left(\sqrt{\frac{j_n 2^{dj_n}}{n}}+2^{-j_n t}\right)+ o(1),\quad a.s.
\end{equation}
\end{theorem}
\textbf{Proof.} In view of decomposition \eqref{dec}, it suffices to handle the first and the last term. The behavior of the second term is given by the previous Proposition \ref{p1}. Let us begin with the first term $R_n(\textbf{u})$. We have for $\textbf{k}=(k_1,\ldots,k_d)\in\mathbb{Z}^d$
\begin{eqnarray*}
\hat{\alpha}_{j_n,\textbf{k}}- \tilde{ \alpha}_{j_n,\textbf{k}}&=&\frac{1}{n}\sum_{i=1}^n \left[\phi_{j_n,\textbf{k}}(F_{1n}(X_{i1}),\ldots, F_{dn}(X_{id}))- \phi_{j_n,\textbf{k}}(F_1(X_{i1}),\ldots, F_d(X_{id}))\right]\\
&=:& \frac{1}{n}\sum_{i=1}^n \xi_{\textbf{k}}(X_{i1},\ldots,X_{id}),
\end{eqnarray*}
where, we set 
\begin{equation}
\xi_{\textbf{k}}(X_{i1},\ldots,X_{id})=\phi_{j_n,\textbf{k}}(F_{1n}(X_{i1}),\ldots, F_{dn}(X_{id}))- \phi_{j_n,\textbf{k}}(F_1(X_{i1}),\ldots, F_d(X_{id})).
\end{equation}
For $d=2$,  \cite{gmt}  observes that,  with $k=(k_1,k_2),$
\begin{equation}\label{fo}
\xi_{k}(X_{i1}, X_{i2})=\xi_{k_1}(X_{i1})\xi_{k_2}(X_{i2})+ \xi_{k_1}(X_{i1})\phi_{j_n k_2}(F_2(X_{i2}))+ \xi_{k_2}(X_{i2})\phi_{j_n k_1}(F_1(X_{i1})), 
\end{equation} 
where $\xi_{k_m}(X_{im})=\phi_{j_n k_m}(F_{mn}(X_{im}))- \phi_{j_n k_m}(F_m(X_{im})),$ for  $m=1,2$.\\
By induction of \eqref{fo}, we obtain for all fixed $d\ge 2$ that
\begin{equation}\label{expr}
\xi_{\textbf{k}}(X_{i1},\ldots,X_{id})=
\sum_{q=0}^{d-1}\sum_{\epsilon_1+\ldots+\epsilon_d=q}\prod_{m=1}^d\left[ \phi_{j_n k_m}(F_m(X_{im}))\right]^{\epsilon_m}\left[ \xi_{k_m}(X_{im})\right]^{1-\epsilon_m},
\end{equation}
where $(\epsilon_1,\ldots\epsilon_d)\in \{0,1\}^d$.
Recall that $\phi_{jl}(u)=2^{j/2}\phi(2^j u- l), \forall j,l\in\mathbb{Z}$. By using the derivability of $\phi$ by hypothesis, we can write, for all $m=1,\ldots,d$
\begin{eqnarray*}
\xi_{k_m}(X_{im})&=& 2^{\frac{j_n}{2}}\phi(2^{j_n} F_{mn}(X_{im})- k_m) - 2^{\frac{j_n}{2}}\phi(2^{j_n} F_{m}(X_{im})- k_m)\\
&=& 2^{\frac{3}{2}j_n} \left[ F_{mn}(X_{im})-F_{m}(X_{im})\right]\phi'(\zeta_{im}),
\end{eqnarray*}
where $\zeta_{im}$ lies between $F_{mn}(X_{im})$ and $F_m(X_{im})$.
Now, combining Chung's (1949) law of the iterated logarithm (LIL)  with
 the boundedness of $\phi$ and $\phi'$, we obtain, for all $m=1,\ldots,d$
\begin{equation}
\vert \xi_{k_m}(X_{im}) \vert\le 2^{\frac{3}{2}j_n}\times \sqrt{\frac{\log\log n}{2n}}\Vert \phi'\Vert_{\infty}, \quad a.s.  
\end{equation}
Thus, for $d=2$, the expression  in \eqref{expr} can be bounded above ; that is 
\begin{equation}\label{fo1}
\vert \xi_{k}(X_{i1},X_{i2}) \vert\le 2^{3j_n}\left( \frac{\log\log n}{2n}\right)\Vert \phi'\Vert^{2}_{\infty}+ 2.2^{2j_n} \sqrt{\frac{\log\log n}{2n}}\Vert \phi'\Vert_{\infty}\Vert \phi\Vert_{\infty}, \quad a.s.  
\end{equation}
Since $$\frac{2^{3j_n}\left( \frac{\log\log n}{2n}\right)}{2^{2j_n} \sqrt{\frac{\log\log n}{2n}}}=\frac{1}{\sqrt{2}}\left( \frac{j_n 2^{2j_n}}{n}\right)^{1/2} \left( \frac{\log\log n}{j_n}\right)^{1/2},$$
which, by hypothesis (H.4) converges to 0, as $n\rightarrow\infty$.
Then  $\;2^{3j_n}\left( \frac{\log\log n}{2n}\right)=o\left(2^{2j_n} \sqrt{\frac{\log\log n}{2n}}\right)$.
That is, for $d=2$, 
$$ \vert \xi_{k}(X_{i1},X_{i2}) \vert= O\left( 2^{2j_n} \sqrt{\frac{\log\log n}{n}}\right).$$
By induction of formula \eqref{fo1},  we have for all $d\ge 2$,
\begin{equation}\label{ineq}
\vert \xi_{\textbf{k}}(X_{i1},\ldots,X_{id}) \vert\le 2^{\frac{3}{2}dj_n}\left( \frac{\log\log n}{n}\right)^{\frac{d}{2}}\Vert \phi'\Vert^{d}_{\infty}+ \cdots+ d.2^{\frac{3}{2}j_n}\sqrt{\frac{\log\log n}{n}} 2^{\frac{d-1}{2}j_n}\Vert \phi'\Vert_{\infty} \Vert \phi\Vert^{d-1}_{\infty}. 
\end{equation}
Note that the number of terms in the summation in the right hand side of inequality \eqref{ineq} is finite. Moreover, as we observe for the case $d=2$, all  these terms  are dominated (\textit{small-o's}) by the last one, which is of order $O\left(2^{\frac{3}{2}j_n}\sqrt{\frac{\log\log n}{n}}2^{\frac{d-1}{2}j_n}\right)$. Then
$$
\vert \xi_{\textbf{k}}( X_{i1},\ldots,X_{id} ) \vert = O\left(2^{\frac{3}{2}j_n}\sqrt{\frac{\log\log n}{n}}2^{\frac{d-1}{2}j_n}\right).
$$
and
\begin{eqnarray*}
|\hat{\alpha}_{j_n,\textbf{k}}- \tilde{ \alpha}_{j_n,\textbf{k}}|&= &\frac{1}{n}\sum_{i=1}^n\vert \xi_{\textbf{k}}(X_{i1},\ldots,X_{id}) \vert \\
&= & O\left(2^{(\frac{2+d}{2})j_n}\sqrt{\frac{\log\log n}{n}} \right)
\end{eqnarray*}
Finally, by using the boundedness of the function $\theta_{\phi}(x)=\sum_{l=1}^{2^{j_n}}|\phi(x-l)|$, we obtain
\begin{eqnarray*}
|\hat{c}_{j_n,\textbf{k}}(\textbf{u})- \tilde{c}_{j_n,\textbf{k}}(\textbf{u})|&\le &\sum_{\textbf{k}\in\{1,\ldots,2^{j_n}\}^d}\vert \hat{\alpha}_{j_n,\textbf{k}}- \tilde{\alpha}_{j_n,\textbf{k}}\vert 2^{\frac{d}{2}j_n}\prod_{m=1}^d\phi(2^{j_n}u_m-k_m) \\
&=& O\left[ \Vert \theta_{\phi}\Vert_{\infty}^d 2^{\frac{2+2d}{2}j_n} \sqrt{\frac{\log\log n}{n}}\right] \\
&=& O\left[ \left(\frac{j_n 2^{(1+d)j_n}}{n}\right)^{1/2} \left(\frac{\log\log n}{j_n}\right)^{1/2}\right] 
\end{eqnarray*}
which, by hypothesis (H.4), converges to 0, as $n\rightarrow\infty$. Hence
\begin{equation}\label{dif}
\sup_{\textbf{u}\in I^d} |R_n(\textbf{u})|\longrightarrow 0, n\rightarrow\infty, \quad a.s.
\end{equation}

To handle the last term $B_n(\textbf{u})$ corresponding to the bias of $\tilde{c}_{j_n}$, we make use of approximation properties in Besov spaces. Let $K_{j_n}$ denote the orthogonal projection kernel onto the sub-space $V_{j_n}$. That is
$$ K_{j_n}(c)(\textbf{u})=\int_{I^d}K_{j_n}(\textbf{u},\textbf{v})c(\textbf{v})d\textbf{v},\quad \textbf{u}\in I^d.$$
Then, we can write
$$B_n(\textbf{u})=\E\tilde{c}_{j_n}(\textbf{u})-c(\textbf{u})= K_{j_n}(c)(\textbf{u})- c(\textbf{u}).$$
 Since  $\phi$ satisfies \textit{Condition} $1(N)$  and $c\in B^t_{\infty,\infty}$, $0<t<N+1$, then applying Theorem 9.4 in \cite{hkt} gives : 
$$\Vert K_{j_n}(c)- c\Vert_{\infty}\le A2^{-j_n t},$$
where $A$ is a positive constant depending on the Besov norm of  $\Vert c\Vert$. Hence
\begin{equation}\label{bias}
\sup_{\textbf{u}\in I^d}\vert B_n(\textbf{u})\vert =O(2^{-j_n t}).
\end{equation}
 Combining  \eqref{dev}, \eqref{dif} and \eqref{bias} gives the proof of the theorem. $\square$\\

\begin{remark}
{\rm Theorem \ref{teo1} implies that if $2^{j_n}\simeq (n/\log n)^{\frac{1}{2t+d}}$, then the rank-based linear estimator $\hat{c}_{j_n}$ achieves the optimal rate of convergence in supremum norm, $(\log n/n)^{\frac{t}{2t+d}}$, over Besov-balls in $B^t_{\infty,\infty}$. This rate is the best possible, as far as the supremum norm loss is concerned $(p=\infty)$, and the  estimated density  is defined on a compact set [see, e.g., \cite{jl} for optimality in minimax theory].  However, notice that this rate is slower than that obtained  in the case of quadratic and pointwise loss functions and established in \cite{gmt} and \cite{aut} for the wavelet linear estimators. }
\end{remark}
\begin{remark}
{\rm For copula densities  in general Besov spaces $ B^{t}_{pq}$ with $t>1/p$, we also have optimal rates for the wavelet linear estimator $\hat{c}_{j_n}$. Indeed, if $c\in B^{t}_{pq}$ with $t>1/p$,  the Sobolev embedding properties entails $B^{t}_{pq}\subset B^{t-1/p}_{\infty,\infty}$. Thus, if $2^{j_n}\simeq (n/\log n)^{\frac{1}{2(t-1/p)+d}}$, then $\hat{c}_{j_n}$ attains the optimal rate : $\E(\Vert \hat{c}_{j_n}-c \Vert_{\infty}) = O\left( (\log n/n)^{\frac{t-1/p}{2(t-1/p)+d}}\right)$ [see, e.g., \cite{dono}, Theorem 1].}
\end{remark}
\begin{remark}
{\rm  As established in \cite{gmt} for the quadratic norm, we have proved by another approach using Chung's (1949) LIL, that the error term associated with the use of ranks (coming from the pseudo-observations) is also negligible for the supremum norm case. That is resorting to pseudo-observations instead of genuine observations does not affect the convergence rate of the linear wavelet estimators of the copula density.}
\end{remark}

\section*{Appendix A : Useful results on empirical process}
\textbf{ Bernstein's inequality (maximal version):} \\
Let $Z_1,\ldots, Z_n$ be independent random variables
with $\mathbb{E}(Z_i)=0, i=1,\ldots,n$ and ${\rm Var}(\sum_{i=1}^{n}Z_i)\le \nu$. Assume further
that for some constant $M >0,\; |Z_i | < M$, $i = 1, . . . , n$. Then for all $t >0$
\begin{equation}
\mathbb{P}\left( \max_{q\le n}\left | \sum_{i=1}^q Z_i \right |  > t\right )\leq 2\exp \left\{\frac{-t^2}{2\nu + (2/3)Mt}\right \},
\end{equation} \\
\bigskip
\textbf{ Lemma A.1 [Einmahl and mason(2000)]:}\\
Let $\mathcal{F}$ and $\mathcal{G}$  be two classes of real-valued measurable
functions on $\mathcal{X}$ satisfying 
$$ \vert f(x)\vert \leq F(x),\quad f\in\mathcal{F},\quad x\in\mathcal{X}, $$
where $F$ is a finite valued measurable envelope function on $\mathcal{X}$;
$$\Vert g\Vert \leq M \quad g\in \mathcal{G},$$
where $M>0$ is a finite constant. Assume that for all probability measure $Q$ with
$0<Q(F^2)<\infty$,
$$ N(\varepsilon(Q(F^2))^{1/2},\mathcal{F},d_Q)\leq C_1 \varepsilon^{-\nu_1}, \qquad 0<\varepsilon<1$$
and 
$$ N(\varepsilon M, \mathcal{G}, d_{Q})\leq C_2 \varepsilon^{-\nu_2}, \qquad 0<\varepsilon<1$$
where $\nu_1, \nu_1, C_1,C_2\geq 1$ are suitable constants. Then we have for all probability measure $Q$ with $0<Q(F^2)<\infty$,
$$ N(\varepsilon M(Q(F^2))^{1/2}, \mathcal{F}\mathcal{G}, d_{Q})\leq C_3 \varepsilon^{-\nu_1-\nu_2}, \qquad 0<\varepsilon<1$$
for some finite constant $0<C_3<\infty$.\\

\textbf{ Proposition 2 [Einmahl and Mason (2000)]:}\\
Let $Z, Z_1,Z_2,\ldots$ , be a sequence of i.i.d. random vectors taking values in $\mathbb{R}^m,m\ge 1$.
For each $n\ge1$, consider the empirical distribution function based on the first $n$ of these random vectors, defined by
$$ G_n(s) = \frac{1}{n}\sum_{i=1}^{n} 1_{Z_i\le s},\quad s\in\mathbb{R}^m,$$
where as usual $z\le s$ means that each component of $z$ is less than or equal to
the corresponding component of $s$. For any measurable real valued function $g$
defined on $\mathbb{R}^m$, set
$$G_n(g) =\int_{\mathbb{R}^m}g(s) dG_n(s), \qquad \mu(g)= \mathbb{E} g(Z)\qquad \text{and} \quad \sigma(g) =Var(g(Z)).$$
Let ${a_n : n\ge 1}$ denote a sequence of positive constants converging to zero. Consider
a sequence $\mathcal{G}_n = \{g^{(n)}_i: i = 1, . . . , k_n\} $ of sets of real-valued measurable
functions on $\mathbb{R}^2$, satisfying, whenever $g^{(n)}_i \in \mathcal{G}_n$:
$$
\mathbb{P}(g^{(n)}_i(Z)=0,\; g^{(n)}_j(Z)) =0,\quad  i\neq j\quad \text{and}\quad \sum_{i=1}^{k_n}\mathbb{P}({g^{(n)}_i (Z) \neq 0})\le 1/2.
$$
Further assume that :\\
For some $0 <r< \infty$, $a_n k_n\rightarrow r$, as $n\rightarrow\infty.$\\
For some $-\infty< \mu_1,\mu_2\infty$, uniformly in $i=1,\ldots,k_n$, for all large $n$, $a_n\mu_1\leq \mu(g^{(n)}_i)\leq a_n\mu_2 $\\
For some $0<\sigma_1<\sigma_2<\infty,$ uniformly in $i = 1,\ldots k_n$, for all large $n$,
$\sigma_1\sqrt{n}a_n \le\bar{\sigma}(g^{(n)}_i)\le \sigma_2\sqrt{n}a_n$\\
For some $0<B<\infty$, uniformly in $i =1,\ldots, k_n$, for all large $n$, $\vert g^{(n)}_i\vert\le B$
\begin{proposition}
	Under these  assumptions, with probability one, for each $0<\varepsilon<1$, there exists $N_{\varepsilon}$ such that for $n\ge N_{\varepsilon}$,
	$$  \max_{1\le i\le k_n}\frac{\sqrt{n} \{G_n(g^{(n)}_i )-\mu(g^{(n)}_i )\}}{\bar{\sigma}(g^{(n)}_i )\sqrt{2|\log a_n|}} \geq 1-\varepsilon.$$
\end{proposition}

\textbf{Talagrand's inequality}:\\
Let $X_i, i=1,\ldots,n$ be an independent and identically distributed  random samples of $X$ with probability law $P$ on $\R$, and $\mathcal{G}$ a $P$-centered (i.e.,$\int gdP=0$ for all $g\in\mathcal{G}$) countable class of real-valued functions on $\R$,
uniformly bounded by the constant {\rm U}. Let $\sigma$ be any positive number such
that $\sigma^2 \ge \sup_{g\in\mathcal{G}}\E(g^2(X))$. Then, Talagrand's (1996) inequality implies that there exists a universal constant $L$ such that for all $t>0$,  
\begin{equation}\label{tal2}
\mathbb{P}\left( \max_{q\le n} \left\Vert \sum_{i=1}^q (g(X_i)  \right\Vert_{\mathcal{G}}  > E+ t \right )\leq L\exp \left\{ \frac{-t}{L{\rm U}}\log (1+\frac{t{\rm U}}{V})\right\},
\end{equation}
with  $$E=\mathbb{E}\left\Vert \sum_{i=1}^n g(X_i) \right\Vert_{\mathcal{G}}\; \text{ and }\; V=\mathbb{E}\left\Vert \sum_{i=1}^n (g(X_i))^2  \right\Vert_{\mathcal{G}}.$$
Further, if $\mathcal{G}$ is a VC-type class of functions, with characteristics $A,\, v$ then, there exist a universal constant $B$ such that : [see, e.g., Gin\'e and Guillou (2001)]
\begin{equation}
E\le B\left[v{\rm U}\log \frac{A{\rm U}}{\sigma}+\sqrt{v}\sqrt{n\sigma^2}\frac{A{\rm U}}{\sigma}\right]
\end{equation} 
Next, if $\sigma<\frac{{\rm U}}{2}$, the constant $A$ may be replaced by 1 at the price of changing the constant $B$, and then if, moreover, $n\sigma^2> C_0\log \left(\frac{{\rm U}}{\sigma}\right)$, we have 
\begin{equation}\label{em}
E\leq C_1\sqrt{n\sigma^2\log \left( \frac{{\rm U}}{\sigma}\right)},\; \text{ and }\; V\le L'n\sigma^2,
\end{equation} where $C_1, L'$ are constants depending only on $A,v,C_0$.
Finally, it follows from \eqref{tal} and \eqref{em} that, for all $t>0$ satisfying : 
 $C_1\sqrt{n\sigma^2 \log \left( \frac{{\rm U}}{\sigma}\right)}\leq t\leq C_2\frac{n\sigma^2}{{\rm U}}$ for all constant $C_2 \ge C_1,$  
\begin{equation}\label{tal1}
\mathbb{P}\left( \max_{n_{k-1}\le n\le n_k} \left\Vert \sum_{i=1}^n g(X_i) \right\Vert_{\mathcal{G}}  > t \right )\leq R\exp \left\{ \frac{-1}{C_3}\frac{t^2}{n\sigma^2}\right\},
\end{equation} where $C_3=\log(1+C_2/L')/RC_2$ and $R$ a constant depending only on $A$ and $v$.

\section*{Appendix B : Proof of Proposition \ref{p1}}
\subsection*{UPPER BOUND}
\begin{lemma}\label{lem1}
	Under the assumptions of Proposition \ref{p1}, one has almost surely
		\begin{equation}
	\limsup_{n\rightarrow\infty} r_n \sup_{\textbf{u}\in I^d}\frac{|D_n(\textbf{u})|}{\sqrt{{ \int_{{\R}^d} K^2(\textbf{x},2^{j_n}\textbf{u})d\textbf{x}}}}\leq \sqrt{\Vert c\Vert_{\infty}}.
\end{equation}
\end{lemma}
\textbf{Proof :}
Given $\lambda>1$, define $n_k=[\lambda^k], k\in\mathbb{N}$, where $[a]$ denotes the integer part of a real $a$. Let  $\delta_m=1/m$, $m\ge 1$ integer, then we can cover the set $I^d$ by a number $l_k$ of small cubes $S_{k,r}$, each of side length $\delta_m 2^{-j_{n_k}}$, with
\begin{equation}\label{lk}
l_k\leq \left(\frac{1}{\delta_m 2^{-j_{n_k}}}+1\right)^d\leq \left(\frac{2}{\delta_m 2^{-j_{n_k}}}\right)^d,
\end{equation}
for $k$ large enough. Let us choose points $\textbf{u}_{k,r}\in S_{k,r}\cap I^d,\; r=1,\ldots,l_k$. We want to prove Lemma \ref{lem1} over the discrete grid of points $\{\textbf{u}_{k,r}: r=1,\ldots,l_k\}$. For all $\eta\in (0,1)$ we claim that
\begin{equation}\label{sss1}
\limsup_{k\rightarrow\infty} \sqrt{\frac{n_k}{(2d\log 2) j_{n_k} 2^{dj_{n_k}}}}\max_{1\le r\le l_k}\max_{n_{k-1}\le n\le n_k }|D_n(\textbf{u}_{k,r})|\le (1+\eta)\sqrt{\Vert c\Vert_{\infty} [K^2]},
\end{equation}
where we note $$ [K^2]= \int_{{\R}^d} \textbf{K}^2(\textbf{x},2^{j_n}\textbf{u})d\textbf{x}.$$
To prove \eqref{sss1}, we apply the maximal version of Bernstein inequality [see, Appendix A above]. Given $\textbf{u}\in I^d$ and $k\in \mathbb{N}$, for all $n$ satisfying : $n_{k-1}\le n\le n_k$ let 
$$ Z_i(\textbf{u})= \textbf{K}(2^{j_{n}}\textbf{U}_i ,2^{j_{n}}\textbf{u})- \mathbb{E}\textbf{K}(2^{j_{n}}\textbf{U}_i ,2^{j_{n}}\textbf{u}),\;\; i=1,\ldots, n. $$
Observe that for each $n$, the $Z_i(\textbf{u})'s$ are independent and identically distributed zero-mean random variables, and for all $\textbf{u}\in I^d$,
\begin{equation}\label{dn}
D_n(\textbf{u})=\frac{2^{dj_n}}{n}\sum_{i=1}^n Z_i(\textbf{u}).
\end{equation}
 By hypothesis (H.2), we have
\begin{eqnarray*}
	\left| \textbf{K}(2^{j_{n}}\textbf{U}_i,2^{j_{n}}\textbf{u}) \right|& =& \prod_{m=1}^d \left|\tilde{K}(2^{j_{n}}U_{im} , 2^{j_{n}}u_m)\right|
	\leq \prod_{m=1}^d \Phi(2^{j_{n}}(U_{im} -u_m)) 
	\leq  \Vert \Phi\Vert^d_{\infty},
\end{eqnarray*}
where $\Vert \Phi\Vert_{\infty}=\sup_{x\in\mathbb{R}}\vert \Phi(x)\vert$.
This implies
$$ \left\vert\mathbb{E}\textbf{K}[(2^{j_{n}}\textbf{U}_i,2^{j_{n}}\textbf{u})]  \right\vert\leq \mathbb{E} \Vert \Phi \Vert^2_{\infty}= \Vert \Phi\Vert^d_{\infty}.$$
Thus, for all  $\textbf{u}\in I^d$,
$$ |Z_i(\textbf{u})|\le 2\Vert \Phi\Vert^d_{\infty}:= M.$$
Since $Z_i(\textbf{u})'s$ are independent and centered, we can write for $n=n_k$
$$Var\left(\sum_{i=1}^{n_k} Z_i(\textbf{u})\right)={n_k} Var(Z_1(\textbf{u}))={n_k}\mathbb{E}(Z_1^2(\textbf{u})).$$
Then using the change of variables $\textbf{s}=2^{-j_{n_k}}\textbf{x}, \textbf{s}=(s_1,\ldots,s_d)$, $\textbf{x}=(x_1,\ldots,x_d)$, we obtain 
\begin{eqnarray*}
	\mathbb{E}(Z_1^2(\textbf{u}))&\leq & \mathbb{E} \textbf{K}^2[(2^{j_{n_k}}\textbf{U}_1,2^{j_{n_k}}\textbf{u})\\
	&\leq & \int_{I^d}\textbf{K}^2(2^{j_{n_k}}\textbf{s},2^{j_{n_k}}\textbf{u})c(\textbf{s})d\textbf{s}\\
	&\leq & 2^{-dj_{n_k}}{\Vert c\Vert_{\infty}}\int_{[0, 2^{j_{n_k}}]^d} \textbf{K}^2(\textbf{x},2^{j_{n_k}}\textbf{u})d\textbf{x}
\end{eqnarray*}
which yields
$$ Var\left(\sum_{i=1}^{n_k} Z_i(\textbf{u})\right)\leq {n_k} 2^{-dj_{n_k}}\Vert c\Vert_{\infty}\int_{\mathbb{R}^d}\textbf{K}^2(\textbf{x},2^{j_n}\textbf{u})d\textbf{x}:=\sigma^2_{k}.$$

Now, applying the maximal version Bernstein's inequality, for each point $\textbf{u}_{k,r}$, we obtain  for all $t>0$,
\begin{equation}
\mathbb{P}\left( \max_{n_{k-1}\le n\le n_k}\left | \sum_{i=1}^n Z_i(\textbf{u}_{k,r}) \right |  > t\right )\leq 2\exp \left\{\frac{-t^2}{2\sigma^2_{k} + (2/3)Mt}\right \},
\end{equation} 
which yields
\begin{eqnarray*}
	\mathbb{P}\left( \max_{1\le r\le l_k} \max_{n_{k-1}\le n\le n_k}\left | \sum_{i=1}^n Z_i(\textbf{u}_{k,r}) \right |  > t\right )&=& \mathbb{P}\left( \bigcup_{r=1}^{l_k} \left\{\max_{n_{k-1}\le n\le n_k}\left | \sum_{i=1}^n Z_i(\textbf{u}_{k,r}) \right |  > t \right \}\right )\\
	&\leq &\sum_{r=1}^{l_k}\mathbb{P}\left( \max_{n_{k-1}\le n\le n_k}\left | \sum_{i=1}^n Z_i(\textbf{u}_{k,r}) \right |  > t\right )\\
	&\leq & l_k 2\exp \left\{ \frac{-t^2}{2\sigma^2_{k} + (2/3)Mt}\right \}.
\end{eqnarray*}
Let $t=\sqrt{2(1+\eta)n_k 2^{-dj_{n_k}}\log 2^{dj_{n_k}}{\Vert c\Vert_{\infty}[ K^2]}}$. 
Then, for $k$ large enough, $t\rightarrow\infty$. Combining this with \eqref{lk},  we obtain after some little algebra,
\begin{eqnarray*}
	\mathbb{P}\left( \max_{1\le r\le l_k} \frac{\max_{n_{k-1}\le n \le n_k }\left | \sum_{i=1}^n Z_i(\textbf{u}_{k,r}) \right |}{\sqrt{2n_k 2^{-dj_{n_k}}\log 2^{dj_{n_k}}{\Vert c\Vert_{\infty}[ K^2]}}}  > \sqrt{1+\eta}\right)
	&\leq & 2 l_k \exp \left\{ \frac{-t^2}{\frac{t^2}{(1+\eta)\log 2^{dj_n}} + \frac{4}{3}\Vert \Phi \Vert^2 t }\right\}\\
	&\leq & 2 l_k \exp \left\{ -(1+\eta)\log 2^{dj_{n_k}} \right\}\\
	&\leq & 2^d \delta^{-d}_m 2^{-d\eta j_{n_k}}.
\end{eqnarray*}
Since the series $\sum_{k\ge 0}2^{-d\eta j_{n_k}}<\infty$, Borel-Cantelli lemma yields
\begin{equation}\label{discr0}
\mathbb{P}\left( \max_{1\le r\le l_k} \frac{ \max_{n_{k-1}\le n \le n_k } \left| \sum_{i=1}^n Z_i(\textbf{u}_{k,r}) \right |}{\sqrt{(2d\log 2)n_k j_{n_k} 2^{-dj_{n_k}}{\Vert c\Vert_{\infty}}[K^2]}}  > \sqrt{1+\eta}\right )
=o(1), 
\end{equation}
That is
\begin{equation}\label{dis}
\limsup_{k\rightarrow\infty}\max_{1\le r\le l_k} \frac{ \max_{n_{k-1}\le n \le n_k } \left| \sum_{i=1}^n Z_i(\textbf{u}_{k,r}) \right|}{\sqrt{(2d\log 2) n_k j_{n_k} 2^{-dj_{n_k}}{\Vert c\Vert_{\infty}}[K^2]}} \le \sqrt{1+\eta},\quad a.s.
\end{equation}
Since the function $x\mapsto x2^{-2x}$ is decreasing for $x>2\log 2$, we have for $n_{k-1}\le n \le n_k $,
and $k$ large enough,
\begin{equation}\label{lam}
\sqrt{\frac{n_k j_{n_k} 2^{-dj_{n_k}}}{n j_{n} 2^{-dj_{n}}}}\le \sqrt{\frac{n_k}{n}}\le \sqrt{\frac{n_k}{n_{k-1}}}\le \sqrt{\lambda}.
\end{equation}
In view of inequality \eqref{lam}, Statement \eqref{dis} yields
\begin{equation}\label{dis1}
\limsup_{k\rightarrow\infty}\max_{1\le r\le l_k} \frac{ \max_{n_{k-1}\le n \le n_k } \left| \sum_{i=1}^n Z_i(\textbf{u}_{k,r}) \right|}{\sqrt{(2d\log 2) n j_{n} 2^{-dj_{n}}{\Vert c\Vert_{\infty}} [K^2]}} \le \sqrt{\lambda(1+\eta)},\quad a.s.
\end{equation}
Now, multiplying  the numerator and the denominator of the fraction in \eqref{dis1} by the factor $2^{dj_{n}}/n$, an recalling the expression of $D_n(\textbf{u})$ in \eqref{dn},  we finally get for all $\eta\in (0,1)$,
\begin{equation}\label{disc}
\limsup_{k\rightarrow\infty} \max_{1\le r\le l_k} \frac{\max_{n_{k-1}\le n \le n_k} \sqrt{n}\left | D_n(\textbf{u}_{k,r}) \right|}{\sqrt{(2d\log 2)j_{n} 2^{dj_{n}}}} \le \sqrt{\lambda(1+\eta)\Vert c\Vert_{\infty} [K^2]},
\end{equation}
which proves Lemma \ref{lem1} over the discrete  grid. \\

Next, to prove Lemma \ref{lem1} between the grid points, we shall make use of  Talagrand's (1996) inequality \eqref{tal2}. 
Let us introduce the sequence of functions  defined as follows : for all $n\ge 1$, $k\ge 1$, $1\le r\le l_k$ and any fixed $\textbf{u}\in S_{k,r}$, define
\begin{equation} \label{gn} g_{k,r}^{(n)}(\textbf{s},\textbf{u})=\textbf{K}(2^{j_{n_k}}\textbf{s},2^{j_{n_k}}\textbf{u}_{k,r}) - \textbf{K}(2^{j_{n}}\textbf{s},2^{j_{n}}\textbf{u}),\qquad \textbf{s}\in I^d. 
\end{equation}
and set, for  all $\lambda>1$,
$$
\mathcal{G}_{k,r}(\lambda)=\left\{g: \textbf{s}\mapsto g_{k,r}^{(n)}(\textbf{s},\textbf{u}):
 \textbf{u}\in S_{k,r}\cap I^d,\; n_{k-1}\le n\le n_k  \right\}.
$$ 
 Let $\textbf{S}=(S_1,\ldots,S_d)$ be a vector  of $[0,1]-$uniform random variables, now we have to check the  following conditions in order to apply Talagrand's (1996) inequality :
\begin{itemize}
	\item [i)] The classes $\mathcal{G}_{k,r}(\lambda), 1\le r\le l_k$, are of VC-type with characteristics $A$ and $v$  ;
	\item[ii)] $\forall g\in\mathcal{G}_{k,r}(\lambda)$, \; $\Vert g\Vert_{\infty} \leq {\rm U} ;$
	\item[iii)] $\forall g\in\mathcal{G}_{k,r}(\lambda)$, \; ${\rm Var}[g(\textbf{S})]\leq \sigma^2_k,$ 
	\item[iv)] $\sigma_k<\frac{{\rm U}}{2}$ and $n_k\sigma^2_k> C_0\log \left(\frac{{\rm U}}{\sigma_k}\right)$, $\;C_0>0$.
\end{itemize}
These conditions will be checked below. \\
Recall that $(\textbf{U}_i=(F_{1i}(X_{i1}),\ldots,F_{di}(X_{id})), i,\ldots,n$ is a sequence of independent and identically distributed vectors of $[0,1]-$uniform components.
We  have shown (see below) that  each class $\mathcal{G}_{k,r}(\lambda)$ satisfies all the conditions i), ii), iii) and iv) for ${\rm U}=2\Vert \Phi\Vert_{\infty}^d$ and $\sigma^2_k=D_0 2^{-dj_{n_k}}{\Vert c\Vert_{\infty}}\omega^2_{\phi}(\delta_m)$, where $\omega_{\phi}$ is the modulus of continuity of $\phi$ defined below \eqref{mc} and $D_0$ is a positive constant depending on $\Vert\Phi\Vert_{\infty}$ and $d$. Then, Talagrand's (1996) inequality  gives : for all $t>0$, 
\begin{equation}\label{tal}
\mathbb{P}\left( \max_{n_{k-1}\le n\le n_k} \left\Vert \sum_{i=1}^n (g(\textbf{U}_i) -\mathbb{E}g(\textbf{U}_i))  \right\Vert_{\mathcal{G}_{k,r}(\lambda)}  >  t \right) \leq R\exp\left\{ \frac{-1}{C_3}\frac{t^2}{n_k\sigma^2_k}\right\},
\end{equation}
which yields, by taking the maximum over $r$ and $t=C_1\sqrt{n_k\sigma^2_k\log \left( \frac{{\rm U}}{\sigma_k}\right)}$,
\begin{equation}
\mathbb{P}\left(\max_{1\le r\le l_k} \max_{n_{k-1}\le n\le n_k} \left\Vert \sum_{i=1}^n (g(\textbf{U}_i) -\mathbb{E}g(\textbf{U}_i))  \right\Vert_{\mathcal{G}_{k,r}(\lambda)}  >  t \right) \leq R l_k\exp\left\{ \frac{-C_1^2}{C_3}\log \left(\frac{{\rm U}}{\sigma_k}\right)\right\}.
\end{equation}
Whenever $m\rightarrow\infty$, $\omega_{\phi}(\delta_m)\rightarrow 0$. Hence, for any $\varepsilon >0$, there exists $m_0\in\mathbb{N}$ such that $\omega_{\phi}(\delta_m)<\varepsilon$ for $m\ge m_0$. Using this fact, we can replace $\sigma_k^2$ by $4D 2^{-dj_{n_k}}\varepsilon{\Vert c\Vert_{\infty}}$, for $m$ large enough.
We also have, for $k$ large enough, 
$$ \log \left(\frac{{\rm U}}{\sigma_k}\right)=\log \left(\frac{{\rm U}}{4D\varepsilon\Vert c\Vert_{\infty}}\right)+j_{n_k}\log 2 \sim j_{n_k}\log 2$$
and thus, for  $k,m$ large enough,
$$ t=C_1\sqrt{n_k\sigma^2_k\log \left( \frac{{\rm U}}{\sigma_k}\right)}\sim\sqrt{4DC_1^2 n_k  2^{-dj_{n_k}}d\, j_{n_k}\log 2\,\varepsilon\Vert c\Vert_{\infty}}.$$
By combining these facts with \eqref{lk} we obtain, with $A_0=\sqrt{2DC_1}$,
\begin{equation}
\mathbb{P}\left(\max_{1\le r\le l_k} \frac{\max_{n_{k-1}\le n\le n_k} \left\Vert \sum_{i=1}^n (g(\textbf{U}_i) -\mathbb{E}g(\textbf{U}_i))\right\Vert_{\mathcal{G}_{k,r}(\lambda)}}{\sqrt{(2d\log 2) n_k j_{n_k} 2^{-dj_{n_k}}}}  > A_0\sqrt{\varepsilon\Vert c\Vert_{\infty}} \right) \leq 4R\delta^{-2} 2^{ -[C_1^2/2C_3-1]2j_{n_k}}.
\end{equation}
Now, we can choose the constant $C_1$ in such a way that $C_1^2/2C_3-1>0$ ; in which case the series $\sum_{k\ge 0} 2^{ -[C_1^2/2C_3-1]2j_{n_k}}$ converges. Thus,  Borel-Cantelli lemma implies 
\begin{equation}\label{bet0}
\mathbb{P}\left(\max_{1\le r\le l_k} \frac{\max_{n_{k-1}\le n\le n_k} \left\Vert \sum_{i=1}^n (g(\textbf{U}_i) -\mathbb{E}g(\textbf{U}_i))\right\Vert_{\mathcal{G}_{k,r}(\lambda)}}{\sqrt{(2d\log 2) n_k j_{n_k} 2^{-dj_{n_k}}}}  > A_0\sqrt{\varepsilon\Vert c\Vert_{\infty}} \right)= o(1),
\end{equation}
that is
\begin{equation}\label{bet1}
\limsup_{k\rightarrow\infty}\max_{1\le r\le l_k} \frac{\max_{n_{k-1}\le n\le n_k} \left\Vert \sum_{i=1}^n (g(\textbf{U}_i) -\mathbb{E}g(\textbf{U}_i))\right\Vert_{\mathcal{G}_{k,r}(\lambda)}}{\sqrt{(2d\log 2) n_k j_{n_k} 2^{-dj_{n_k}}}} \le A_0\sqrt{\varepsilon\Vert c\Vert_{\infty}},\quad a.s.
\end{equation}

Arguing as in the discrete case, with Statement \eqref{lam} in view,  we  conclude that
\begin{equation}\label{bet1}
\limsup_{k\rightarrow\infty}\max_{1\le r\le l_k} \frac{\max_{n_{k-1}\le n\le n_k} \sqrt{n}\left\Vert D_n(\textbf{u}_{k,r}) -D_n(\textbf{u}))\right\Vert_{\mathcal{G}_{k,r}(\lambda)}}{\sqrt{(2d\log 2)  j_{n} 2^{dj_{n}}}}\le A_0\sqrt{\lambda\varepsilon\Vert c\Vert_{\infty}},\quad a.s.,
\end{equation}
which completes the proof of Lemma \eqref{lem1} between the grid-points.\\

Now recapitulating, we can infer from \eqref{disc} and \eqref{bet1} that
\begin{equation}\label{discr}
\limsup_{k\rightarrow\infty} \max_{1\le r\le l_k} \frac{\max_{n_{k-1}\le n \le n_k} \sqrt{n}\left | D_n(\textbf{u}) \right|}{\sqrt{(4\log 2)j_{n} 2^{-2j_{n}}}} \le \sqrt{\lambda(1+\eta)\Vert c\Vert_{\infty}[K^2]}+A_0\sqrt{\lambda\varepsilon\Vert c\Vert_{\infty}},\quad a.s.
\end{equation}
Since $\eta$ and $\varepsilon$ are arbitrary, letting $\lambda\rightarrow 1$ completes the proof of  Lemma \ref{lem1}.$\square$ 

\section*{ Checking conditions i), ii), iii), iv)}
\textbf{Checking  i):} Observe that the elements of the class $\mathcal{G}_{k,r}(\lambda)$ may be rewritten as
$$
g_{k,r}^{(n)}(\textbf{s},\textbf{u})=
\prod_{m=1}^d \widetilde{K}(2^{j_{n_k}}s_m,2^{j_{n_k}}u_{k,r,m})- \prod_{m=1}^d \widetilde{K}(2^{j_{n}}s_m,2^{j_{n}}u_{m}),
$$
where $\widetilde{K}(x,y)=\sum_{l=1}^{2^{j_n}} \phi(x-l)\phi(y-l), $ with $\phi$ compactly supported and is of bounded variation. For $m=1,\ldots,d$, define the classes of functions :
$ \mathcal{F}_m=\{v\mapsto \sum_{l\in\mathbb{Z}} \phi(2^{j}w-l)\phi(2^{j}v-l): w\in [0,1], j\in\mathbb{N} \}$.
By Lemma 2 in \cite{gn}, $\mathcal{F}_1,\ldots,\mathcal{F}_m$ are VC-type classes of functions. Moreover, $ \mathcal{F}_1,\ldots,\mathcal{F}_m$ are uniformly bounded.
Indeed,  for all $w\in[0,1],\; j\in \mathbb{N}$, we have
$
\left\vert \sum_{l=1}^{2^{j_n}} \phi(2^{j}w-l)\phi(2^{j}\cdot-l)\right\vert \leq \Vert \phi \Vert_{\infty} \Vert \theta_{\phi}\Vert_{\infty}$, as the function  $\theta_{\phi}(x)=\sum_{l=1}^{2^{j_n}} |\phi(x-l)|$ is bounded.   By Lemma A.1 in 
\cite{em}, this implies that the product $\mathcal{F}_1 \cdots\mathcal{F}_m$ is also a VC-type class of functions.  Now, using properties (iv) and (v) of Lemma 2.6.18 in 
\cite{vdv}, we can infer that the classes of functions $\mathcal{G}_{k,r}(\lambda)$ are of VC-type for all $k,r$ fixed.\\

\textbf{Checking  ii):} For all $k\ge 1,\; 0\le r\le l_k,\; n_{k-1}\le n\le n_k$,  using hypothesis (H.2), we can write 
\begin{eqnarray*}
	\left\vert g_{{k},r}^{(n)}(\cdot,\textbf{u}) \right\vert & \leq & \left|\textbf{K}( 2^{j_{n_k}}\cdot,2^{j_{n_k}}\textbf{u}_{{k},r})\right| + \left|\textbf{K}( 2^{j_n}\cdot,2^{j_n}\textbf{u})\right| \\
	&\leq & \prod_{m=1}^d \widetilde{K}(2^{j_{n_k}}\cdot ,2^{j_{n_k}}u_{{k},r,m}) + \prod_{m=1}^d \widetilde{K}(2^{j_{n}}\cdot,2^{j_{n}}u_{m}) \leq  2\Vert \Phi \Vert^d_{\infty}
\end{eqnarray*}
and ii) holds with ${\rm U} =2\Vert \Phi \Vert^d_{\infty}$.\\

\textbf{Checking  iii):} For all $k\ge 1,\; 0\le r\le l_k,\; n_{k-1}\le n\le n_k$. As in \cite{gn} we choose $\lambda\in (0,1),$ such that $j_{n_k}=j_n$. By a change of variable $\textbf{s}=\textbf{u}+2^{-j_{n_k}}\textbf{x},\; \textbf{s}=(s_1,\ldots,s_d),\; \textbf{x}=(x_1,\ldots,x_d),\;\textbf{u}=(u_1,\ldots,u_d)$, we have 
\begin{eqnarray*}
\mathbb{E}\left[ \left(g_{{k},r}^{(n)}(\textbf{S},\textbf{u})\right)^2\right] & = &  \mathbb{E}\left[ \left(\textbf{K}(2^{j_{n_k}}\textbf{S}, 2^{j_{n_k}}\textbf{u}_{k,r}) - \textbf{K}(2^{j_{n_k}}\textbf{S}, 2^{j_{n_k}}\textbf{u})\right)^2\right]   \\
	& = & \int_{I^d} \left(\textbf{K}(2^{j_{n_k}}\textbf{s},2^{j_{n_k}}\textbf{u}_{k,r})- \textbf{K}(2^{j_{n_k}}\textbf{s},2^{j_{n_k}}\textbf{u})\right)^2 c(\textbf{s})d\textbf{s}  \\
	& \le &  \frac{\Vert c\Vert_{\infty}}{2^{-dj_{n_k}}}\int_{[-2^{j_{n_k}},2^{j_{n_k}}]^d} \left(\textbf{K}(2^{j_{n_k}}\textbf{u}+\textbf{x}, 2^{j_{n_k}}\textbf{u}_{k,r}) - \textbf{K}(2^{j_{n_k}}\textbf{u}+\textbf{x}, 2^{j_{n_k}}\textbf{u})\right)^2 d\textbf{x}. \\
& \le &  {2^{-dj_{n_k}}}{\Vert c\Vert_{\infty}}\int_{\mathbb{R}^d} \left(\textbf{K}(2^{j_{n_k}}\textbf{u}+\textbf{x}, 2^{j_{n_k}}\textbf{u}_{k,r}) - \textbf{K}(2^{j_{n_k}}\textbf{u}+\textbf{x}, 2^{j_{n_k}}\textbf{u})\right)^2 d\textbf{x}. 
\end{eqnarray*} 
To simplify, let us take $\textbf{w}=2^{j_{n_k}}\textbf{u}+\textbf{x}$, then $d\textbf{x}=d\textbf{w}$, $w=(w_1,\ldots,w_d)\in\mathbb{R}^d$. \\
Put $A(\textbf{w})=\textbf{K}(\textbf{w},2^{j_{n_k}}\textbf{u}_{k,r}) - \textbf{K}(\textbf{w},2^{j_{n_k}}\textbf{u})$  ; using the multiplicativity of the kernel $\textbf{K}$, we can rewrite $A(\textbf{w})$ as 
\begin{eqnarray*}
A(\textbf{w})&= & \prod_{l=1}^d\widetilde{K}(w_l,2^{j_{n_k}}u_{k,r,l})- \prod_{l=1}^d\widetilde{K}(w_l,2^{j_{n_k}}u_{k,r,l})\\
&=& \sum_{l=1}^d \left[\widetilde{K}(w_l,2^{j_{n_k}}u_{k,r,l})-\widetilde{K}(w_l,2^{j_{n_k}}u_{l})\right]\prod_{p=1,p\neq l}^d\widetilde{K}(w_p,2^{j_{n_k}}u_{k,r,p})
\end{eqnarray*}
For any $\delta>0$, the modulus of continuity of $\phi$ is defined as
\begin{equation}\label{mc} \omega_{\phi}(\delta)=\{ \sup |\phi(x) -\phi(y)|: |x-y|\le \delta\}. \end{equation}
Recall that  
$ \widetilde{K}(x,y)=\sum_{h=1}^{2^{j_n}} \phi(x-h)\phi(y-h).$ 
Combining these facts with the inequality $(a_1+\cdots+a_d)^2\le d (a_1^2+\cdots +a_d^2)$, and  Fubini's Theorem, we get 
\begin{eqnarray*}
	\int_{\mathbb{R}^d}\vert A(\textbf{w})\vert^2 d\textbf{w} &\leq& d\int_{\mathbb{R}^d}\sum_{l=1}^d \left[\widetilde{K}(w_l, 2^{j_{n_k}}u_{k,r,l})-\widetilde{K}(w_l,2^{j_{n_k}}u_{l})\right]^2 \prod_{p=1,p\neq l}^d\widetilde{K}^2(w_l,2^{j_{n_k}}u_{k,r,p})d\textbf{w}.
\end{eqnarray*}	
Then
\begin{eqnarray*}
	\int_{\mathbb{R}^d}\vert A(\textbf{w})\vert^2 d\textbf{w} &\leq & d\sum_{l=1}^d \int_{\mathbb{R}^d}\left[\sum_{h=1}^{2^{j_n}}\phi(2^{j_{n_k}}w_l-h)[\phi(2^{j_{n_k}}u_{k,r,l}-h) - \phi(2^{j_{n_k}}u_l-h)] \right]^2\\ 
	\hspace{3cm} & &\qquad\qquad \times\prod_{p=1,p\neq l}^d\widetilde{K}^2(w_p,2^{j_{n_k}}u_{k,r,p}) d\textbf{w}\\
	&\le &d\omega_{\phi}^2(\delta_m)\sum_{l=1}^d \int_{\mathbb{R}^d}\left[\sum_{h}^{2^{j_n}}\phi(2^{j_{n_k}}w_l-h)\right]^2\prod_{p=1,p\neq l}^d\widetilde{K}^2(w_p,2^{j_{n_k}}u_{k,r,p}) d\textbf{w}\\
	&\le &d\omega_{\phi}^2(\delta_m)\sum_{l=1}^d \int_{\mathbb{R}}\left[\sum_{h=1}^{2^{j_n}}\phi(2^{j_{n_k}}w_l-h)\right]^2dw_l \prod_{p=1,p\neq l}^d \int_{\mathbb{R}}\widetilde{K}^2(w_p,2^{j_{n_k}}u_{k,r,p}) d{w_p}.
\end{eqnarray*}
 Now, since the family $\{\phi(\cdot-h): h=1,\ldots,{2^{j_n}}\}$ is an orthonormal basis, the quantity $\int_{\mathbb{R}}\left(\sum_{h=1}^{2^{j_n}}\phi(w_l-h)\right)^2dw_l$ can be bounded by a constant $M_0$ ; thus
\begin{eqnarray*}
\int_{\mathbb{R}^d}\vert A(\textbf{w})\vert^2 d\textbf{w} &\le & M_0 d\omega_{\phi}^2(\delta_m)\sum_{l=1}^d  \prod_{p=1,p\neq l}^d \int_{\mathbb{R}}\widetilde{K}^2(w_p,2^{j_{n_k}}u_{k,r,p}) d{w_p}\\
	&\leq & M_0 d^2\omega^2_{\phi}(\delta_m)D,
\end{eqnarray*}
where we use Hypothesis (H.2) for the last inequality, with $D$ a positive constant depending on $\Vert\Phi\Vert_{\infty}$. 
Finally, we obtain
\begin{equation}
\mathbb{E}\left[ \left(g_{{k},r}^{(n)}(\textbf{S},\textbf{u})\right)^2\right] \leq M_0 d^2 D 2^{-dj_{n_k}}{\Vert c\Vert_{\infty}}\omega^2_{\phi}(\delta_m), 
\end{equation}
and iii) holds with \begin{equation} \label{ss1}\sigma^2_k=D_0 2^{-dj_{n_k}}{\Vert c\Vert_{\infty}}\omega^2_{\phi}(\delta_m), \qquad D_0=M_0 d^2 D.\\
\end{equation}

\textbf{Checking  iv):} For $m>0$ fixed, we have 
$$\frac{\sigma_k}{\rm U}= \frac{D_0^{1/2}2^{-\frac{d}{2}j_{n_k}}{\Vert c\Vert_{\infty}^{1/2}}\omega_{\phi}(\delta_m)}{2\Vert \Phi \Vert^d_{\infty}}\rightarrow 0, k\rightarrow\infty,$$
which implies that $\frac{\sigma_k}{\rm U}<\varepsilon$, for all $\varepsilon>0$ and $k$ large enough. Hence, for $\varepsilon=1/2$, we have $\sigma_k<\frac{\rm U}{2}$. We also have, for all large $k$,
$$ \frac{n_k\sigma_k^2}{\log \left( \frac{\rm U}{\sigma_k}\right)}= \frac{D_0 n_k{\Vert c\Vert_{\infty}}\omega^2_{\phi}(\delta_m)}{j_{n_k}2^{d j_{n_k}}\,\log 2}\longrightarrow\infty,$$
by  Hypothesis (H.4). This readily implies that, for any constant $C_0>0$,
$n_k\sigma^2_k >C_0\log \left( \frac{\rm U}{\sigma_k}\right)$
for all large $k$, and iv) holds. $\square$

 \subsection*{LOWER BOUND}
 \begin{lemma}\label{lem2}
	Under the assumptions of Proposition \ref{p1}, one has almost surely
		\begin{equation}\label{ss3}
	\liminf_{n\rightarrow\infty} r_n \sup_{\textbf{u}\in I^d}\frac{|D_n(\textbf{u})|}{\sqrt{{ \int_{\mathbb{R}^d} K^2[(\textbf{x},2^{j_n}\textbf{u})]d\textbf{x}}}}\geq \sqrt{\Vert c\Vert_{\infty}}.
\end{equation}
\end{lemma}
\textbf{Proof :} It is an adaptation of the proof of Proposition 2 in  \cite{gn}, which is, itself, inspired by Proposition 2 in \cite{em}. According to this latter proposition, \eqref{ss3} holds if and only if for all $\tau >0$, and all large $n$, there exists $k_n=:k_n(\tau)$ points $\textbf{z}_{i,n}=(z_{1,i,n},\ldots, z_{d,i,n})\in I^d, \; i=1,\ldots,k_n$ such that, for functions
$ g_{i}^{(n)} (\textbf{s})=\textbf{K}(2^{j_n}\textbf{s}, 2^{j_n}\textbf{z}_{i,n}),\; \textbf{s}\in I^d,$ and for $\textbf{U}=(U_1,\ldots,U_d)$ a random vector with joint density $c$, the following conditions hold :
\begin{itemize}
	\item[C.1)] $ \mathbb{P}(g_{i}^{(n)}(\textbf{U})\neq 0,\; g_{i'}^{(n)}(\textbf{U})\neq 0)=0,\quad \forall i\neq i' ;$ 
	\item[C.2)] $\sum_{i=1}^{k_n}\mathbb{P}(g_{i}^{(n)}(\textbf{U})\neq 0)\leq 1/2 ;$ 
	\item[C.3)] $2^{-j_n}k_n\longrightarrow r\in ]0,\infty[ ;$ 
	\item[C.4)] $\exists\; \mu_1, \mu_2\in \mathbb{R}$ : $2^{-dj_n}\mu_1\leq \mathbb{E}g_{i}^{(n)}(\textbf{U})\leq 2^{-dj_n}\mu_2, \quad \forall i=1,\ldots,k_n ;$ 
	\item[C.5)] $\exists\; \sigma_1, \sigma_2>0$ : $2^{-dj_n}\sigma_1^2\leq {\rm Var}[g_{i}^{(n)}(\textbf{U})]\leq 2^{-dj_n}\sigma_2^2, \quad \forall i=1,\ldots,k_n ;$ 
	\item[C.6)] $ \Vert g_{i}^{(n)}\Vert_{\infty}<\infty,\quad \forall i=1,\ldots,k_n ; \; \forall n\geq 1;$ 
\end{itemize}

Now, we have to check these conditions. By hypothesis the copula density $c$ is continuous and bounded on $I^d$, then there exists some orthrotope $D\subset I^d$ such that
$ \max_{\textbf{s}\in D} c(\textbf{s})= \Vert c\Vert_{\infty}$.
Thus, for all $\tau >0$ there exists $\textbf{s}_0\in D$ such that $c(\textbf{s}_0)\geq (1-\tau)\Vert c\Vert_{\infty}$. Let
\begin{equation}\label{d} D_{\tau}=\{\textbf{s}\in D : c(\textbf{s})\geq (1-\tau)\Vert c\Vert_{\infty}\}, \end{equation}
and choose a subset $ D_0 \subset D_{\tau}\; \text { such that } \; \mathbb{P}(\textbf{U}\in D_0)\leq \frac{1}{2}$.  Suppose that
$ D_0=\prod_{j=1}^d[a_j,b_j],$ with $ 0\le a_j<b_j\le 1, \;  \; \text{ and } b_j-a_j=\ell,\;\forall j=1,\ldots,d. $\\
Set $\delta=3B$ and define
$$ z_{j,i,n}=a+i\delta 2^{-j_n}, \quad i=1,\ldots,\left[ \frac{b-a}{\delta 2^{-j_n}} \right]-1:=k_n,\quad j=1,\ldots,d$$
where $[x]$ designs the integer part of a real $x$.

\textbf{Checking C.1) :} Recall that $\phi$ is supported on $[0,B]$, then
$$ 
g_{i}^{(n)}(\textbf{U})\neq 0 \iff \forall k,l\in \mathbb{Z} \quad \left\{ \begin{array}{ccc} 0 \le & 2^{j_n}U_j-l & \le B,\; j=1,\ldots,d \quad (1) 
\\0 \le & 2^{j_n}z_{j,i,n}-l &\le B, \; j=1,\ldots,d \quad (2) \end{array}\right.$$
and
$$ 
g_{i'}^{(n)}(\textbf{U})\neq 0 \iff \forall k,l\in \mathbb{Z}\quad  \left\{ \begin{array}{ccc} 0 \le & 2^{j_n}U_j-l & \le B, \; j=1,\ldots,d \quad(1)' \\0 \le & 2^{j_n}z_{j,i',n}-l &\le B, \; j=1,\ldots,d\quad (2)'\end{array}\right. $$
Combining (2) and $(2)'$ gives, for every $j=1,\ldots,d$, 
$$|z_{j,i,n} -z_{j,i',n}|\leq 2^{-j_n}B.\quad (3)$$ But, by definition, for all $i\neq i'$, 
$|z_{j,i,n} -z_{j,i',n}|> \delta 2^{-j_n}=3B2^{-j_n}$, which contradicts (3). Hence, the event $ \{g_{i}^{(n)}(\textbf{U})\neq 0,g_{i'}^{(n)}(\textbf{U})\neq 0\}$ is empty for $i\neq i'$ and condition C.1) holds.\\

\textbf{Checking C.2) :} For all $n\geq 1$, the sets $\{g_{i}^{(n)}(\textbf{U})\neq 0\},\, i=1,\ldots, k_n$ are disjoint in view of Condition C.1). Then, we have
$$\sum_{i=1}^{k_n}\mathbb{P}(\{g_{i}^{(n)}(\textbf{U})\neq 0\})=\mathbb{P}\left(\bigcup_{i=1}^{k_n}\{g_{i}^{(n)}(\textbf{U})\neq 0\}\right).$$
Now, it suffices to show that
$ \bigcup_{i=1}^{k_n}\{g_{i}^{(n)}(\textbf{U})\neq 0\}\subset\{ \textbf{U}\in D_0 \}$.
From statements (1) and (3) above, we can write, fro all $j=1,\ldots,d$
\begin{eqnarray*}
	&-B\leq 2^{j_n}(U_j-u_{j,i,n}) \leq B & \\
	& u_{j,i,n} -2^{-j_n}B\leq U_j \leq u_{j,i,n}+ 2^{-j_n}B & \\
	&a_j\le a_j+ (3i-1)2^{-j_n}B\leq U_j \leq a_j+ (3i+1)2^{-j_n}B \le b_j.&
\end{eqnarray*}
That is $U_j\in [a_j,b_j]$, and hence $\textbf{U}=(U_1,\ldots,U_d)\in \prod_{j=1}^d [a_j,b_j]=D_0.$ It follows  that,  
\begin{eqnarray*}
& \forall\;i=1,\ldots,k_n,\quad \{g_{i}^{(n)}(\textbf{U})\neq 0\} \subset \{\textbf{U}\in D_0\}&\\
& \bigcup_{i=1}^{k_n}\{g_{i}^{(n)}(\textbf{U})\neq 0 \}\subset \{\textbf{U}\in D_0\}&\\
& \mathbb{P}\left( \bigcup_{i=1}^{k_n}\{g_{i}^{(n)}(\textbf{U})\neq 0 \}\right) \leq \mathbb{P}(\{\textbf{U}\in D_0\})\leq \frac{1}{2}.&
\end{eqnarray*}
Hence, C.2) is fulfilled.\\

\textbf{Checking C.3):} It is immediate, since
$$ 2^{-j_n}k_n=2^{-j_n}\left( \left[ \frac{b-a}{\delta 2^{-j_n}} \right]-1\right)=\left[ \frac{b-a}{\delta} \right]-2^{-j_n}\rightarrow \left[ \frac{b-a}{\delta} \right]=:r>0,\; n\rightarrow\infty.$$

\textbf{Checking C.4) :} Using a change of variables $\textbf{s}=2^{-j_n}\textbf{x},\; \textbf{s}=(s_1,\ldots,s_d),\; \textbf{x}=(x_1,\ldots,x_d)$, we have
\begin{eqnarray*}
	\vert \mathbb{E}g_{i}^{(n)}(\textbf{U})\vert &\leq & \int_{I^d_{\epsilon}}\left\vert \textbf{K}(2^{j_n}\textbf{s}, 2^{j_n}\textbf{z}_{i,n})\right\vert c(\textbf{s})d\textbf{s}\\
	&\leq &2^{-dj_n} \Vert c\Vert_{\infty} \int_{\mathbb{R}^d}\left\vert \textbf{K}(\textbf{x}, 2^{j_n}\textbf{z}_{i,n})\right\vert d\textbf{x}\\
	&\leq &2^{-dj_n} \Vert c\Vert_{\infty} \int_{\mathbb{R}^d}\vert\prod_{j=1}^d \widetilde{K}(x_j,2^{j_n}u_{j,i,n})\vert  dx_j\\
	&\leq &2^{-dj_n} \Vert c\Vert_{\infty} \int_{\mathbb{R}^d}\prod_{j=1}^d\Phi(x_j-2^{j_n}u_{j,i,n})dx_j\\
	&\leq & 2^{-2j_n}\mu,
\end{eqnarray*}
where $\mu= \Vert c\Vert_{\infty}\int_{\mathbb{R}^d}\prod_{j=1}^d\Phi(x_j-2^{j_n}u_{j,i,n})dx_j$ exists, because the function $\Phi$ is integrable by hypothesis (H.2). 
The last inequality is equivalent to $$ -2^{-2j_n}\mu\leq \mathbb{E}g_{i}^{(n)}(U,V) \leq  2^{-2j_n}\mu,\; \forall i=1,\cdots,k_n. $$
That is C.4) holds.\\

\textbf{Checking C.5) :} For $n\geq 1,\; i=1,\ldots,k_n$, using a change of variables $\textbf{s}=2^{-j_n}\textbf{x}+\textbf{z}_{i,n},\; \textbf{s}=(s_1,\ldots,s_d),\; \textbf{x}=(x_1,\ldots,x_d),\; \textbf{z}_{i,n}=(z_{1,i,n},\ldots,z_{d,i,n})$, we can write
\begin{eqnarray*}
	{\rm Var}[g_{i}^{(n)}(\textbf{U})] &\leq & \mathbb{E}\left[\left(g_{i}^{(n)}(\textbf{U})\right)^2\right]\\
	&\leq & \int_{I^d} \textbf{K}^2(2^{j_n}\textbf{s}, 2^{j_n}\textbf{z}_{i,n}) c(\textbf{s})d\textbf{s}\\
	&\leq & 2^{-dj_n}{\Vert c\Vert_{\infty}}\int_{\mathbb{R}^d} \textbf{K}^2(\textbf{x}+2^{j_n}\textbf{z}_{i,n}, 2^{j_n}\textbf{z}_{i,n}) d\textbf{x}.
\end{eqnarray*}
Putting
$ \sigma^2_2:={\Vert c\Vert_{\epsilon}}\int_{\mathbb{R}^d} \textbf{K}^2(\textbf{x}+2^{j_n}\textbf{z}_{i,n},2^{j_n}\textbf{z}_{i,n}) d\textbf{x}$ yields $${\rm Var}[g_{i}^{(n)}(\textbf{U})]\leq 2^{-dj_n}\sigma^2_2,$$ which is the upper bound in condition C.5).
For the lower bound, we have 
\begin{eqnarray*}
	& &{ \rm Var}[g_{i}^{(n)}(\textbf{U})] = \mathbb{E}\left[\left(g_{i}^{(n)}(\textbf{U})\right)^2\right] -\left[\mathbb{E}g_{i}^{(n)}(\textbf{U})\right]^2\\
	&=& \int_{I^d} \textbf{K}^2(2^{j_n}\textbf{s}, 2^{j_n}\textbf{z}_{i,n}) c(\textbf{s})d\textbf{s} - \left( \int_{I^d} \textbf{K}(2^{j_n}\textbf{s}, 2^{j_n}\textbf{z}_{i,n}) c(\textbf{s})d\textbf{s}\right)^2.
\end{eqnarray*}
Put $\mu_{n}^2=\left( \int_{I^d} \textbf{K}(2^{j_n}\textbf{s}, 2^{j_n}\textbf{z}_{i,n}) c(\textbf{s})d\textbf{s}\right)^2.$ 
Noting that  $D_{\tau}\subset I^d$, by a change of variables $\textbf{x}=2^{j_n}\textbf{s}$, we obtain
\begin{eqnarray*}
	{\rm Var}[g_{i}^{(n)}(\textbf{U})] &\geq & \int_{D_{\tau}} \textbf{K}^2(2^{j_n}\textbf{s}, 2^{j_n}\textbf{z}_{i,n}) c(\textbf{s})d\textbf{s} - \mu_{n}^2\\
	&\geq &(1-\tau)\Vert c\Vert_{\infty}\int_{D_{\tau}} \textbf{K}^2(2^{j_n}\textbf{s}, 2^{j_n}\textbf{z}_{i,n}) d\textbf{s} - \mu_{n}^2\\
	&\geq &(1-\tau)\Vert c\Vert_{\infty}2^{-dj_n}\int_{\mathbb{R}^d} \textbf{K}^2(\textbf{x}, 2^{j_n}\textbf{z}_{i,n})d\textbf{x} - \mu_{n}^2.
\end{eqnarray*}
Proceeding again to the same change of variables, and observing  from hypothesis (H.3) that \\
$\int_{\mathbb{R}^d} \textbf{K}(\textbf{x}, 2^{j_n}\textbf{z}_{i,n}) d\textbf{x}=1$, we can write
\begin{eqnarray*}
	\mu_{n}^2 &\le & \left(\Vert c\Vert_{\infty}2^{-dj_n} \int_{\mathbb{R}^d} \textbf{K}(\textbf{x}, 2^{j_n}\textbf{z}_{i,n}) d\textbf{x}\right)^2\le \Vert c\Vert^2_{\infty}2^{-2dj_n},
\end{eqnarray*}
which implies $-\mu_n^2\geq -\Vert c\Vert^2_{\epsilon}2^{-4j_n}.$ Thus, for $n$ large enough, we obtain the lower bound in condition C.5), i.e.
\begin{eqnarray*} 
	{\rm Var}[g_{i}^{(n)}(\textbf{U})] &\geq & 2^{-dj_n}(1-\tau)\Vert c\Vert_{\infty}\int_{\mathbb{R}^d} \textbf{K}^2(\textbf{x}, 2^{j_n}\textbf{z}_{i,n})d\textbf{x}- \Vert c\Vert_{\infty}^2 2^{-2dj_n}\\
	&\geq & 2^{-dj_n} \sigma_1^2 + o(1),
\end{eqnarray*}
with $\sigma_1^2:=(1-\tau)\Vert c\Vert_{\infty}\int_{\mathbb{R}^d} K^2(\textbf{x}, 2^{j_n}\textbf{z}_{i,n}]d\textbf{x}.$ Finally, C.5) holds.\\
Moreover, letting $\tau\rightarrow 0$, we get 
$\sigma_2^2=\sigma_1^2=\Vert c\Vert_{\infty}\int_{\mathbb{R}^d} \textbf{K}^2(\textbf{x}, 2^{j_n}\textbf{u}d\textbf{x}$.\\

\textbf{Checking C.6) :} For all $\textbf{s}\in I^d$,$n\geq 1,\quad i=1,\ldots, k_n$, by using  hypotheses (H.1-2) and the multiplicativity of kernel $\textbf{K}$, we have
\begin{eqnarray*}
\vert  g_{i}^{(n)}(\textbf{s})\vert &=& \vert\textbf{K}(2^{j_n}\textbf{s}, 2^{j_n}\textbf{z}_{i,n})\vert= \prod_{m=1}^d\vert\widetilde{K}(2^{j_n}s_m, 2^{j_n}z_{i,n,m})\vert \\
 &\leq &\prod_{m=1}^d \sum_{l=1}^{2^{j_n}}\vert \phi(2^{j_n}s_m-l)\phi(2^{j_n}z_{i,n,m}-l)\vert \\
	&\leq & \Vert\phi \Vert_{\infty}^d \prod_{m=1}^d\sum_{l=1}^{2^{j_n}}\vert \phi(2^{j_n}s_m-l)\vert \\ 
	&\leq &\Vert\phi \Vert_{\infty}^d\Vert\theta_{\phi} \Vert_{\infty}^d.
\end{eqnarray*}
Hence, $\sup_{n\ge 1, 1\le i\le k_n}\Vert g_{i}^{(n)}\Vert \leq \Vert\phi \Vert_{\infty}^d\Vert\theta_{\phi} \Vert_{\infty}^d $, and C.6) holds.\\
Since  Conditions C.1-2-3-4-5-6) are fulfilled, we can now apply Proposition 2 in \cite{em} to complete the proof of Lemma \ref{lem2}. $\square$\\

Finally, Lemma \ref{lem1} and Lemma \ref{lem2} give the proof of Proposition \eqref{p1}.$\square$
\section*{Conflict of interest}
 On behalf of all authors, the corresponding author states that there is no conflict of interest.\\
 

%
%



\end{document}